\documentclass{tran-l}
\usepackage{amscd,amssymb,amsmath,verbatim,hyperref,
epsfig,latexsym,subfigure}
\usepackage[arrow,matrix,graph,frame,poly,arc,tips]{xy}

\newtheorem{thm}{Theorem}[section]
\newtheorem{cor}[thm]{Corollary}
\newtheorem{lemma}[thm]{Lemma}
\newtheorem{prop}[thm]{Proposition}

\theoremstyle{definition}
\newtheorem{definition}[thm]{Definition}
\newtheorem{example}[thm]{Example}
\newtheorem{remark}[thm]{Remark}

\newtheorem{conj}{Conjecture}

\numberwithin{equation}{section}

\newcommand{\abs}[1]{\left|#1\right|}
\newcommand{\surj}{\twoheadrightarrow}
\newcommand{\inj}{\hookrightarrow}

\newcommand{\A}{{\mathcal A}}

\newcommand{\RR}{{\mathcal R}}

\newcommand{\B}{{\mathfrak B}}
\newcommand{\I}{{\mathfrak I}}
\newcommand{\m}{{\mathfrak m}}
\newcommand{\p}{{\mathfrak p}}
\newcommand{\q}{{\mathfrak q}}

\newcommand{\Z}{{\mathbb Z}}
\newcommand{\Q}{{\mathbb Q}}
\newcommand{\R}{{\mathbb R}}
\newcommand{\C}{{\mathbb C}}
\renewcommand{\P}{{\mathbb P}}
\renewcommand{\k}{\Bbbk}

\newcommand{\PP}{{\operatorname{Ass}(\RR^1(\A))}}

\DeclareMathOperator{\Tor}{Tor}
\DeclareMathOperator{\Ext}{Ext}
\DeclareMathOperator{\Hom}{Hom}
\DeclareMathOperator{\Sym}{Sym}
\DeclareMathOperator{\Hilb}{Hilb}
\DeclareMathOperator{\HP}{P}

\DeclareMathOperator{\rank}{rank}

\DeclareMathOperator{\gr}{gr}

\DeclareMathOperator{\im}{im}
\DeclareMathOperator{\coker}{coker}
\DeclareMathOperator{\codim}{codim}
\DeclareMathOperator{\spn}{span}
\DeclareMathOperator{\ann}{ann}
\DeclareMathOperator{\lin}{lin}
\DeclareMathOperator{\id}{id}
\DeclareMathOperator{\reg}{reg}

\DeclareMathOperator{\ch}{char}
\DeclareMathOperator{\cx}{cx}

\begin{document}

\title[Resonance, syzygies, Chen groups, and BGG]%
{Resonance, linear syzygies, Chen groups, and the
Bernstein-Gelfand-Gelfand correspondence}

\author[Henry K. Schenck]{Henry K. Schenck}
\address{Department of Mathematics,
Texas A\&M University, College Station, TX 77843}
\email{\href{mailto:schenck@math.tamu.edu}{schenck@math.tamu.edu}}
\urladdr{\href{http://www.math.tamu.edu/~schenck/}%
{http://www.math.tamu.edu/\~{}schenck}}

\author[Alexander I. Suciu]{Alexander~I.~Suciu}
\address{Department of Mathematics,
Northeastern University,
Boston, MA 02115}
\email{\href{mailto:a.suciu@neu.edu}{a.suciu@neu.edu}}
\urladdr{\href{http://www.math.neu.edu/~suciu/}%
{http://www.math.neu.edu/\~{}suciu}}

\thanks{Both authors were supported by NSF 
Collaborative Research grant DMS 03-11142; 
Schenck was also supported by NSA grant MDA 904-03-1-0006 
and ATP grant 010366-0103.}



\begin{abstract}
If $\A$ is a complex hyperplane arrangement, with complement $X$,
we show that the Chen ranks of $G=\pi_1(X)$ are equal to the
graded Betti numbers of the linear strand in a minimal, free resolution
of the cohomology ring $A=H^*(X,\k)$, viewed as a module
over the exterior algebra $E$ on $\A$:
\[
\theta_k(G) = \dim_{\k}\Tor^E_{k-1}(A,\k)_k, \quad \text{for $k\ge 2$},
\]
where $\k$ is a field of characteristic $0$.  The 
{\it Chen ranks conjecture} asserts that, for $k$ sufficiently large, 
$\theta_k(G) =(k-1) \sum_{r\ge 1} h_r \binom{r+k-1}{k}$, 
where $h_r$ is the number of $r$-dimensional 
components of the projective resonance variety $\RR^{1}(\A)$.  
Our earlier work on the resolution of $A$ over $E$ and the 
above equality yield a proof of the conjecture for graphic 
arrangements. Using results on the geometry of $\RR^{1}(\A)$ and a 
localization argument, we establish the inequality
\[
\theta_k(G) \ge (k-1) \sum_{r\ge 1} h_r  \binom{r+k-1}{k}, 
\quad \text{for $k\gg 0$}, 
\]
for arbitrary $\A$.  
Finally, we show that there is a polynomial $P(t)$ of 
degree equal to the dimension of  $\RR^1(\A)$, such that 
$\theta_k(G) = P(k)$, for all $k\gg 0$. 
\end{abstract}
\maketitle

\section{Introduction}
\label{intro}

\subsection{Orlik-Solomon algebra}
\label{subs:coho}
Let $\A=\{H_1,\dots ,H_n\}$ be an arrangement of complex hyperplanes in $\C^{\ell}$.  
A fundamental question in the 
subject is to decide whether a 
given topological invariant of the complement, 
$X(\A)=\C^{\ell}\setminus \bigcup_{H\in \A}\, H$, is determined by
the intersection lattice, $L(\A)=\{\bigcap_{H\in \A'}\, H \mid
\A'\subseteq \A\}$, and, if so, to find an explicit combinatorial 
formula for such an invariant. 

For example, in \cite{OS}, 
Orlik and Solomon showed that the cohomology ring of the complement 
is entirely determined by 
$L(\A)$.  More precisely, the Orlik-Solomon algebra
$A=H^*(X(\A),\Z)$ is the quotient of the exterior algebra
$E=\bigwedge (\Z^n)$ on generators $e_1, \dots , e_n$
in degree $1$ by the ideal $I$ generated by all elements of
the form $\partial e_{i_1\dots i_r}:=\sum_{q}(-1)^{q-1}e_{i_1} \cdots
\widehat{e_{i_q}}\cdots e_{i_r}$, for which
$\codim H_{i_1}\cap \cdots \cap H_{i_r} < r$.  
Notice that $I$ is generated in degrees $2$ and higher; 
in particular, $A_0=E_0=\Z$ and $A_1=E_1=\Z^n$.  

For each element $a\in A_1$, the Orlik-Solomon algebra can 
be turned into a cochain complex $(A,a)$.  The $i^{\text{th}}$ 
term of this complex is simply the degree $i$ graded piece of $A$, 
and the differential is given by multiplication by $a$: 
\begin{equation}
\label{aomoto}
(A,a)\colon \quad 
\xymatrix{
0 \ar[r] &A_0 \ar[r]^{a} & A_1
\ar[r]^{a}  & A_2 \ar[r]^{a}& \cdots \ar[r]^{a} 
& A_{\ell}\ar[r] & 0}.
\end{equation}
This complex arose in the work of Aomoto \cite{Ao} 
on hypergeometric functions, and in the work of 
Esnault, Schechtman and Viehweg \cite{ESV} 
on cohomology with coefficients in local systems.  
In \cite{Yuz}, Yuzvinsky showed that, for generic $a$, 
the complex \eqref{aomoto} is exact.    

\subsection{Resonance varieties}
\label{subs:res}
In \cite{Fa}, Falk initiated the study of the cohomology 
jumping loci for an arrangement complement. 
Fix a field $\k$; abusing notation, we will also 
denote by $A=H^*(X(\A),\k)$ the Orlik-Solomon 
algebra over $\k$.  
The {\em resonance varieties} of $\A$ are the loci of points 
$a=\sum_{i=1}^na_ie_i \leftrightarrow (a_1:\dots :a_n)$ in 
$\P(A_1) \cong \P^{n-1}$ for which $(A,a)$ fails to be exact.  
More precisely, for each $k\ge 1$,
\[
\RR^{k}(\A)=\{a \in \P^{n-1} \mid  H^k(A,a)\ne 0\}.
\]

The resonance varieties of $\A$ lie in the hyperplane
$\sum_{i=1}^n a_i = 0$, and depend only on the 
lattice-isomorphism type of $L(\A)$. In \cite{Fa}, Falk 
also introduced the notion of a {\em neighborly
partition}: a partition $\Pi$ of $\A$ is neighborly if, for
any rank two flat $Y\in L_2(\A)$ and any block $\pi$ of $\Pi$,
\begin{equation}
\label{eq:np}
\mu(Y) \le \abs{Y\cap \pi} \Longrightarrow Y\subseteq \pi, 
\end{equation}
where $\mu$ denotes the M\"obius function of the lattice.  
Now assume $\ch\k=0$.  Falk showed that all components 
of $\RR^1(\A)$ arise from neighborly partitions of 
sub-arrangements of $\A$. In particular, each flat 
$Y \in L_2(\A)$ gives rise to a ``local" component 
of dimension $\mu(Y)-1$; see Example \ref{K4} 
for how this works. 

Falk also conjectured that the components of $\RR^{1}(\A)$ 
are projective linear subspaces; this was proved in \cite{CScv}, 
and generalized to $\RR^{k}(\A)$ in \cite{CO}.  
Libgober and Yuzvinsky \cite{LY} showed that the components 
of $\RR^{1}(\A)$ are in fact disjoint, and positive-dimensional.  
These facts will be used in an essential way later in the paper. 
We note that the characteristic zero assumption is necessary; 
the aforementioned results depend on this.  A thorough treatment of 
resonance varieties over arbitrary fields (and even commutative 
rings) can be found in a recent paper of Falk \cite{Fa04} (see also 
\cite{MS}).  

\subsection{Lower central series and Chen ranks}
\label{subs:lcs chen}
The fundamental group of the complement, $G(\A)=\pi_1(X(\A))$, 
is not necessarily determined by the intersection lattice.  Even so, 
the ranks of the lower central series quotients,
\[
\phi_k(G):=\rank G_k/G_{k+1}, 
\]
are combinatorially determined (here $G_1=G$, $G_2=G'=[G,G]$, 
and $G_{k+1}=[G,G_{k}]$).  In fact, due to the formality of $X(\A)$, 
the LCS ranks depend only on the $\k$-algebra 
$A=H^*(X(\A),\k)$, where $\ch\k=0$.  
Explicit formulas for the LCS ranks of an arrangement 
group are available in some cases 
(most notably, when $L(\A)$ is supersolvable), but a general, 
all-encompassing formula remains elusive.  
See \cite{FR}, \cite{Su}, \cite{Yuz1}  for surveys of the problem, 
and \cite{SS}, \cite{PS1} for recent developments.  
 
More manageable topological invariants are the Chen ranks, 
\[
\theta_k(G):=\phi_k(G/G'').
\] 
Introduced by K.T.~Chen in his thesis \cite{Ch51}, these are 
the LCS ranks of the maximal metabelian quotient of $G$. 
For example, if $F_n$ is the free group of rank $n$, 
then $\theta_1(F_n)=n$, and 
$\theta_k(F_n)=(k-1)\binom{k+n-2}{k}$, for $k\ge 2$. 
The original motivation for the study of the invariants 
$\theta_k(G)$ came from the study of link groups, see 
Massey's paper \cite{Ma} and references therein. 

The study of Chen ranks of arrangement groups was 
started in \cite{CS1}, \cite{CSai}. The ranks $\theta_k(G(\A))$ 
were determined in a number of cases, including the pure 
braid groups.   They proved to be quite subtle and useful 
invariants, distinguishing in some instances 
between arrangement groups with same LCS ranks.  

\subsection{Resonance and Chen ranks}
\label{subs:reschen}
For about a decade, it was an open question whether 
the Chen ranks of an arrangement are combinatorially 
determined, see e.g.~\cite[\S2.3]{FR}. This question was 
recently settled in the affirmative in \cite{PS}.  
There remained the question of computing explicitly 
the Chen ranks of an arrangement group $G(\A)$ 
in terms of the intersection lattice $L(\A)$.  Based 
on the work in \cite{CS1}, \cite{CSai}, a precise 
combinatorial formula for the Chen ranks 
was conjectured in \cite{Su}:

\begin{conj}[Resonance formula for Chen ranks]
\label{conj:chen}
Let $G=G(\A)$ be an arrangement group, and let 
$h_r$ be the number of components
of $\RR^{1}(\A)$ of dimension $r$.  Then, for $k \gg 0$:
\begin{equation}
\label{eq:chenconj}
\theta_k(G)= (k-1) \sum_{r\ge 1} h_r  \binom{r+k-1}{k}.  
\end{equation}
\end{conj}

In other words,  
$\theta_k(G)= \sum_{r\ge 1} h_r  \theta_k(F_{r+1})$.  
This formula can easily be verified for a pencil 
of $n$ lines (with $G=F_{n-1}\times \Z$), a near-pencil 
(with $G=F_{n-2}\times \Z^2$), or a product of 
such arrangements. Much less obviously, the 
conjecture holds for the braid arrangements \cite{CS1}, and for  
``decomposable" arrangements  \cite{CSai}, \cite{PS1}.   
Our methods give a unified proof of Conjecture \ref{conj:chen} 
for all these classes of arrangements, and in fact apply 
more generally.

It was originally conjectured in \cite{Su} that 
equality \eqref{eq:chenconj} holds for all $k\ge 4$, 
but Example \ref{ex:ceva3} below shows 
this is false.  It turns out that the value for which 
$\theta_k(G)$ is given by a fixed polynomial in $k$ 
depends on the Castelnuovo-Mumford regularity 
of the linearized Alexander invariant. 

\subsection{Resonance and linear syzygies}
\label{subs:ressyz}
In \cite{SS} we observed that there was a close connection 
between the resonance variety $\RR^{1}(\A)$ and 
the linear syzygies of $A$, where the Orlik-Solomon algebra 
$A$ is viewed as a module over the exterior algebra $E$.  
We define the {\em linear strand} in a minimal free resolution of
$A$ over $E$ as the subcomplex of the form
\begin{equation*}
\label{Ares}
\xymatrixcolsep{14pt}
\xymatrix{
0& A \ar[l]  & E \ar[l] & E^{\beta_{12}}(-2) \ar[l] 
& E^{\beta_{23}} (-3)\ar[l]& \cdots \ar[l] 
&E^{\beta_{i,i+1}} (-i-1) \ar[l]& \cdots \ar[l]
}
\end{equation*}
where $\beta_{ij}=\dim_{\k} \Tor^E_i(A,\k)_j$. While $\beta_{12}$
is determined by the M\"obius function of $L(\A)$, a purely 
combinatorial formula for $\beta_{i, i+1}$ is unknown.  
However, many examples suggested to us that: 

\begin{conj}[Resonance formula for the linear strand]
\label{conj:linconj}
For $k\gg 0$, the graded Betti numbers of the linear strand
are given by
\begin{equation}
\label{eq:linconj}
\beta_{k-1,k}= (k-1) \sum_{r\ge 1} h_r  \binom{r+k-1}{k} .
\end{equation}
\end{conj}

\subsection{Outline and Results}
\label{subs:results}
We now outline the  structure of paper,  and state our 
main results.  

We start in \S\ref{cast} with a review of the linearized Alexander 
invariant $\B$.    First considered in \cite{CSai}, this graded module 
over the symmetric algebra $S = \Sym(A_1^{*})$  is closely related 
to both the resonance variety $\RR^1(\A)$ \cite{CScv}, and to the 
Chen ranks $\theta_k(G)$ \cite{PS}.  

In \S\ref{sect:BGG}, we express the linearized Alexander invariant 
as an $\Ext$ module.  Key to our approach is the 
Bernstein-Gelfand-Gelfand correspondence, which gives 
(for our purposes) a relationship between linear complexes 
over the exterior algebra $E$ (on generators $e_1,\dots , e_n$) 
and graded modules over the symmetric algebra $S$ 
(on generators $x_1,\dots , x_n$).  In Theorem \ref{bext2}, we show: 
\begin{equation}
\B \cong\Ext_S^{\ell-1}(F(A), S).
\end{equation}
Here $F(A)$ is the top cohomology module of the Aomoto complex
$(A,a)$ tensored with $S$, with $a=\sum_{i=1}^n x_ie_i$, 
cf.~\eqref{fa}. The module $F(A)$ was first studied by 
Eisenbud-Popescu-Yuzvinsky 
in \cite{EPY}.  The formulation of the BGG correspondence 
in \cite{EFS} shows that the local cohomology modules
of $F(A)$ determine the free resolution of $A$ over $E$.

As a first application, we prove in  Corollary \ref{thetator} that 
Conjecture \ref{conj:chen} and Conjecture \ref{conj:linconj} 
are equivalent:
\begin{equation}
\label{eq:thetabeta}
\theta_k(G) = \beta_{k-1,k},\quad  \text{for all $k\ge 2$}.  
\end{equation}

As another application, we prove in Theorem \ref{graphic} 
that the Chen ranks conjecture holds for graphic arrangements. 
If $\Gamma$ is a graph and $\A=\A({\Gamma})$ the corresponding 
arrangement, we show that:
\begin{equation}
\label{eq:thetagraph} 
\theta_k(G) = (k-1)(\kappa_2 + \kappa_3), \quad  \text{for all $k\ge 3$},  
\end{equation}
where $\kappa_s$ is the number of cliques of size $s+1$.  On the 
other hand, the components of $\RR^1(\A)$ are all $2$-dimensional, 
and there is one component for each triangle or complete quadrangle 
in $\Gamma$.   Thus, formula \eqref{eq:thetagraph} agrees 
with \eqref{eq:chenconj}. 

In \S \ref{commalg} we use the results of \cite{CScv}, \cite{LY}, \cite{PS}  
and some commutative algebra to prove that the Chen ranks have 
polynomial growth, controlled by the resonance variety: there exists 
a polynomial $P(t) \in \mathbb{Q}[t]$, of degree 
equal to the dimension of $\RR^{1}(\A)$, such that 
$\theta_k(G) = P(k)$, for all $k \gg 0$. In particular, this implies  
\begin{equation}
\label{eq:ratio} 
\lim_{k\to\infty} \frac{\theta_k(G)}{(k-1) 
\sum_{r\ge 1}  h_r  \binom{r+k-1}{k}} \in \mathbb{Q}.
\end{equation}
As an easy corollary of this, we compute the complexity of 
the Orlik-Solomon algebra $A$, viewed as a module over $E$, 
in the case when $\ell=3$. 

In \S \ref{sect:exseq}, we use a localization argument to give 
a (sharp) lower bound on the Chen ranks, thereby  
proving one direction of Conjecture \ref{conj:chen}.  
For an arbitrary arrangement group $G$, we show (in 
Corollary \ref{halfchen}):
\begin{equation}
\label{eq:thetaineq}
\theta_k(G)\ge (k-1) \sum_{r\ge 1} h_r  \binom{r+k-1}{k}, 
\quad \text{for all $k \gg 0$}.
\end{equation}

In \S \ref{sect:examples}, we give some examples illustrating 
the fact that the Chen ranks formula does not hold for small values 
of $k$.  This phenomenon can be interpreted in terms of 
certain local cohomology modules, which reflect subtle 
combinatorial behavior in $L(\A)$.  

\section{Dramatis Personae}
\label{cast}
Let $\A$ be an arrangement of complex hyperplanes in $\C^{\ell}$, 
with complement $X(\A)$.
Since we are primarily interested in the fundamental group
$G(\A)=\pi_1(X(\A))$, we may restrict our attention to affine line arrangements
in $\C^2$.  Indeed, if $\A$ is an arbitrary arrangement, let $\A'$ 
be a generic two-dimensional section of $\A$.  Then, by the Lefschetz-type 
theorem of Hamm and L\^{e} \cite{HL}, the inclusion $X(\A')\inj X(\A)$ induces 
an isomorphism $\pi_1(X(\A'))\xrightarrow{\simeq} \pi_1(X(\A))$. 

Now, if $\A$ is an (affine) line arrangement in $\C^2$, the cone on $\A$ 
is a central plane arrangement in $\C^3$, with complement 
$X(\A)\times (\C\setminus 0)$, and fundamental group $G(\A)\times \Z$, 
see \cite{OT}.  In view of the above, we will assume throughout 
that $\A$ is a central arrangement.  Most often, this will be an 
arrangement in $\C^3$, with projectivization a line arrangement in $\P^2$.

\subsection{Chen groups and Alexander invariant}
\label{subsect:chenalex}
Let $G'=[G,G]$ be the derived subgroup of $G$, and $G''=(G')'$ 
the second derived subgroup. The group $G/G'$ is the maximal abelian 
quotient of $G$, whereas $G/G''$ is its maximal metabelian quotient.
The {\em $k$-th Chen group} of $G$ is, by definition, the $k$-th
lower central series quotient of $G/G''$.
Let $\theta_k(G)=\rank \gr_k (G/G'')$ be its rank.
For example, if $G=F_n$, the free group of rank $n$, then, 
as shown by Chen \cite{Ch51}:
\begin{equation}
\label{eq:chenfn}
\theta_k(F_n)=(k-1)\cdot\binom{k+n-2}{k} , \quad \text{for $k\ge 2$}.
\end{equation}

Fix a basis $\{t_1,\dots ,t_n\}$ for
$H_1(X)=G/G'=\Z^n$, corresponding to meridians around
the hyperplanes in $\A=\{H_1,\dots, H_n\}$.  This identifies the 
group ring $\Z[G/G']$ with the ring of Laurent polynomials
$\Lambda =\Z[t_1^{\pm 1},\dots, t_n^{\pm 1}]$.
Let $\I$ be the augmentation ideal of $\Lambda$, and let
$\gr \Lambda$ be the graded ring associated to the $\I$-adic filtration.
The map $t_i\mapsto 1-x_i$ identifies this ring with the polynomial
ring $S=\Z[x_1,\dots , x_n]$.

The {\em Alexander invariant}, $B=B(\A)$, is the abelian group $G'/G''$, 
endowed with the $\Lambda$-module structure induced by the conjugation 
action in the extension $0\to G'/G'' \to G/G'' \to G/G' \to 0$.   
As shown in \cite{CSai}, the module $B$ admits a finite presentation
of the form 
\begin{equation}
\label{eq:Bpres}
\xymatrix{
(A_2 \oplus E_3)\otimes \Lambda\ar[r]^(.6){\Delta}
& E_2\otimes \Lambda \ar[r] & B \ar[r] & 0
}.
\end{equation}

Let $\gr B=\bigoplus _{k\ge 0}\I^{k}B/\I^{k+1}B$ 
be the associated graded module. 
A basic observation of W.S. Massey \cite{Ma} asserts that
$\gr_k (G/G'')= \gr_{k-2} B$, for $k\ge 2$. Consequently,
\begin{equation}
\label{eq:massey}
\sum_{k\ge 0} \theta_{k+2}\, t^k = \Hilb(\gr B\otimes \k), 
\end{equation}
for any field $\k$ of characteristic $0$.

\subsection{Linearized Alexander invariant}
\label{subs:linalex} 
In general, the module $B$ is hard to compute, depending as it 
does on finding a (braid monodromy) presentation for the group $G$,
and carrying out a laborious Fox calculus algorithm.  The resulting
presentation matrix $\Delta$ typically involves very complicated
Laurent polynomials.  On a conceptual level, it is not at all clear
whether $B$ is combinatorially  determined, since $G$ is not always
determined by $L(\A)$.
For all these reasons, it is convenient to look at a simplified
version of the Alexander invariant, which carries
all the essential information we want to extract from this module.

The {\em linearized Alexander invariant}, $\B=\B(\A)$, is the
graded $S$-module presented by the ``linearization" of the 
matrix $\Delta$:
\begin{equation}
\label{eq:newlinB}
\xymatrix{
(A_2 \oplus E_3)\otimes S\ar[r]^(.6){\Delta^{\lin}}
& E_2\otimes S \ar[r] & \B \ar[r] & 0
}.
\end{equation}
The matrix $\Delta^{\lin}$ appears in the statement
of Theorem 4.6 and in Remark 4.7 from \cite{CScv}; see also 
\cite{MS} and \cite{PS} for more general contexts.  
The reason for the terminology is as follows:  
Viewing $\Lambda \otimes \C $ as the coordinate ring of the algebraic 
torus $(\C^*)^n$ and $S \otimes \C$ as the coordinate ring of $\C^n$, 
the entries of $\Delta^{\lin}$ are the derivatives at 
$\mathbf{1}\in (\C^*)^n$ of the corresponding entries of $\Delta$.  

Let us describe the linearized Alexander matrix in a concrete fashion, 
following \cite{CScv}. If $\delta_i$ denotes the $i$-th differential in 
the Koszul complex, then
\begin{equation}
\label{eq:presB}
\Delta^{\lin}=\alpha_2\otimes \id +\delta_3,
\end{equation}
where $\alpha_2$ is the adjoint of the canonical projection
$\gamma_2\colon E_2\surj A_2$, given by 
$\alpha_2(Y,j)=e_j \wedge \sum_{i\in Y} e_i$, 
where $\{(Y,j)\mid Y \in  L_2(\A) \text{ and } j\in Y\setminus \min Y\}$
is the usual ${\bf nbc}$-basis for $A_2$, see e.g.~\cite{OT}.

\subsection{The module $\B$, resonance, and Chen ranks}
\label{subsec:ann}
When $\k$ is a field of characteristic zero (which will be our standing 
assumption from now on), the linearized Alexander invariant 
$\B:=\B \otimes \k$ is related in a very concrete way to 
both the resonance variety $\RR^1(\A)$, and to the 
fundamental group $G=G(\A)$.  

The first important fact about the module $\B$ 
is the identification of the variety defined by its annihilator ideal.  

\begin{thm}[\cite{CScv}]
\label{fact1}
The annihilator of the linearized Alexander invariant defines
the resonance variety:
\begin{equation}
\label{eq:vann}
V(\ann \B)=\RR^1(\A).
\end{equation}
\end{thm}

The second important fact about the module $\B$ is the
following linearized version of Massey's result. 

\begin{thm}[\cite{PS}] 
\label{fact2}
The Chen ranks, $\theta_k=\theta_k(G)$, $k\ge 2$, are equal to the 
dimensions of the graded pieces of the linearized Alexander invariant:
\begin{equation}
\label{eq:linmas}
\sum_{k\ge 2} \theta_{k}\, t^k = \Hilb(\B).
\end{equation}
In particular, the Chen ranks are combinatorially determined, 
and the Hilbert polynomial $\HP(\B,t) \in \mathbb{Q}[t]$ 
gives the asymptotic Chen ranks: for $k \gg 0$, $\theta_k =\HP(\B,k)$. 
\end{thm}
\begin{example} [Braid arrangement]
\label{K4}

Let $\A$ be the braid arrangement in $\P^2$, with defining 
polynomial $Q=xyz(x-y)(x-z)(y-z)$.  
 From the matroid (see Figure~\ref{fig:k4}), 
it is easy to see that the Orlik-Solomon algebra
$A$ is the quotient of the exterior algebra $E$ on generators 
$e_0,\dots,e_5$ by 
the ideal $I = \langle \partial e_{145}, \partial e_{235}, 
\partial e_{034}, \partial e_{012}, \partial e_{ijkl} \rangle$, 
where $ijkl$ runs over all four-tuples; it turns out that the elements 
$\partial e_{ijkl} $ are redundant. 

The (minimal) free resolution of $A$ as a module over $E$ begins:
\begin{equation*}
\label{braidres}
\xymatrixcolsep{18pt}
\xymatrix{
0& A \ar[l]  & E \ar[l] & E^{4}(-2) \ar[l]_(.55){\partial_1} 
& E^{10} (-3)\ar[l]_(.45){\partial_2} &E^{15}(-4) \oplus E^6(-5)  
\ar[l]_(.6){\partial_3} & \cdots \ar[l]
}, 
\end{equation*}
where $\partial_1=\begin{pmatrix}
\partial e_{145}&\partial e_{235}&
\partial e_{034}&\partial e_{012}\end{pmatrix}$, 
and $\partial_2$ is equal to
\[
{\Small
\begin{pmatrix}
\xymatrixrowsep{1pt}
\xymatrixcolsep{1pt}
\xymatrix{
e_1-e_4 & e_1-e_5 & 0 & 0 & 0 & 0 & 0 & 0 & e_3-e_0 & e_2-e_0\\
0 & 0 & e_2-e_3 & e_2-e_5 & 0 & 0 & 0 & 0 & e_0-e_1 & e_0-e_4\\
0 & 0 & 0 & 0 & e_0-e_3 & e_0-e_4 & 0 & 0 & e_1-e_5 & e_2-e_5\\
0 & 0 & 0 & 0 & 0 & 0 & e_0-e_1 & e_0-e_2 & e_3-e_5 & e_4-e_5
}
\end{pmatrix}\!.
}
\]
\begin{figure}
\subfigure{%
\label{fig:k4arr}%
\begin{minipage}[t]{0.3\textwidth}
\setlength{\unitlength}{14pt}
\begin{picture}(5,3.5)(-3,-0.3)
\multiput(0,1)(0,2){2}{\line(1,0){4}}
\multiput(1,0)(2,0){2}{\line(0,1){4}}
\put(0,4){\line(1,-1){4}}
\put(0,0){\line(1,1){4}}
\put(3,-0.5){\makebox(0,0){$0$}}
\put(4.5,-0.25){\makebox(0,0){$1$}}
\put(4.5,1){\makebox(0,0){$2$}}
\put(4.5,3){\makebox(0,0){$4$}}
\put(4.5,4.3){\makebox(0,0){$3$}}
\put(1,4.5){\makebox(0,0){$5$}}
\end{picture}
\end{minipage}
}
\subfigure{%
\label{fig:k4mat}%
\setlength{\unitlength}{0.6cm}
\begin{minipage}[t]{0.4\textwidth}
\begin{picture}(5,3.5)(-2,-0.2)
\put(3,3){\line(1,-1){3}}
\put(3,3){\line(-1,-1){3}}
\put(1.5,1.5){\line(3,-1){4.5}}
\put(4.5,1.5){\line(-3,-1){4.5}}
\multiput(0,0)(6,0){2}{\circle*{0.3}}
\multiput(1.5,1.5)(3,0){2}{\circle*{0.3}}
\multiput(3,3)(0,-2){2}{\circle*{0.3}}
\put(0,0){\makebox(-1,0){$5$}}
\put(1,1.5){\makebox(0,0.5){$4$}}
\put(3,3.5){\makebox(0,0){$1$}}
\put(3,1){\makebox(0,-1){$3$}}
\put(5,1.5){\makebox(0,0.5){$2$}}
\put(6,0){\makebox(1,0){$0$}}
\end{picture}
\end{minipage}
}
\caption{\textsf{The braid arrangement and its matroid}}
\label{fig:k4}
\end{figure}

The resonance variety $\RR^1(\A)\subset \P^5$ has $4$ 
local components, corresponding to the triple points, and 
$1$ essential component (i.e., one that does not come from 
any proper sub-arrangement), corresponding to the neighborly 
partition $\Pi=(05|13|24)$:
\[
\begin{aligned}
&  \{ x_1 + x_4 + x_5=x_0=x_2=x_3=0 \} ,\ 
\{ x_2 + x_3 + x_5=x_0=x_1=x_4=0 \} ,  \\
&\{ x_0 + x_3 + x_4= x_1=x_2=x_4=0 \}, \
 \{ x_0 + x_1 + x_2=x_3=x_4=x_5=0 \} ,\\
& \{ x_0+x_1+x_2=x_0- x_5=x_1-x_3=x_2-x_4=0 \} .
\end{aligned}
\]

The linearized Alexander invariant 
$\B=\coker\left(  \Delta^{\lin}\colon S^{31} \to S^{15} \right) $ 
is a module over the ring $S=\k[x_0,\dots,x_5]$; 
the presentation matrix $ \Delta^{\lin}$ can be reduced by row and 
column operations to the matrix $\vartheta\colon S^{14} \to S^4$, 
with transpose 
\[
\begin{pmatrix}
\sum x_i&0&0&0\\
0&\sum x_i&0&0\\
0&0&\sum x_i&0\\
0&0&0&\sum x_i\\
x_1 + x_4 + x_5&0&0&0\\
0&x_2 + x_3 + x_5&0&0\\
0&0&x_0 + x_3 + x_4&0\\
0&0&0&x_0 + x_1 + x_2\\
-x_3&0&0&x_3\\
0&x_4&0&x_4\\
0&0&-x_5&x_5\\
x_0&x_0&0&0\\
x_2&0&-x_2&0\\
0&x_1&x_1&0
\end{pmatrix}.
\]

A computation shows that 
\[
\Hilb(\B)= (4t^2+2t^3-t^4)/(1-t)^2 =
4t^2+10t^3+15t^4+20 t^5  +\cdots.
\]
Thus,  $\theta_1 = 6$,  $\theta_2 = 4$, 
and $\theta_k = 5(k-1)$, for $k \ge 3$. 
This agrees with the computations in \cite{CS1}, and 
with the values predicted by Conjecture ~\ref{conj:chen}. 
A generalization will be given in Theorem~\ref{graphic}. 
\end{example}

\section{The Bernstein-Gelfand-Gelfand correspondence 
and $H^*(A,a)$}
\label{sect:BGG}

\subsection{The BGG correspondence}
\label{subs:bgg}
Let $V$ be a finite-dimensional vector space over a field $\k$. In this 
section, we connect our cast of characters. The key tool is the
Bernstein-Gelfand-Gelfand correspondence, which is an isomorphism
between the category of linear free complexes over the exterior
algebra $E=\bigwedge (V)$ and the category of graded free 
modules over the symmetric algebra $S=\Sym(V^*)$.
An introduction to the BGG correspondence may be found in 
Chapter 7 of \cite{e03}; additional sources are \cite{EFS}, \cite{DE},
and \cite{EPY}. 

Notice that if we take $\Sym(V^*)$ to be generated in degree 
one, then $\bigwedge (V)$ is generated in degree $-1$. The 
convention in arrangement theory (and the convention of this paper)
is that the exterior algebra is generated in degree $1$. To
distinguish between gradings, we write $E'$ for an exterior algebra
with generators in degree $-1$, and $E$ if the generators are in
degree $1$. 

Let $\mathbf{L}$ denote the functor from the category of
graded $E'$-modules to the category of linear free complexes
over $S$, defined as follows: for a graded $E'$-module $P$, 
$\mathbf{L}(P)$ is the complex
\begin{equation}
\label{eq:lp}
\xymatrix{
\cdots \ar[r]& P_i \otimes S \ar[r]^(.45){d_i} &  P_{i-1} \otimes S
 \ar[r]& \cdots
}
\end{equation}
with differentials $d_i\colon p\otimes 1 \rightarrow \sum e_ip \otimes x_i$. 
We will be applying the functor $\mathbf{L}$ to the $E'$-module 
$A'  = E'/I \otimes E'(-\ell)$, where $I$ is the Orlik-Solomon ideal; 
tensoring with $E'(-\ell)$ shifts so the unit of the algebra is in degree
$\ell$, and the generators are in degree $\ell-1$, with $\ell$ the
dimension of the ambient space of the arrangement. 

Similarly, let $\mathbf{R}$ denote the 
functor from the category of graded $S$-modules to the 
category of linear free complexes over $E'$: for a graded $S$-module
$M$, $\mathbf{R}(M)$ is
the complex
\begin{equation}
\label{eq:rp}
\xymatrix{
\cdots \ar[r]& \Hom_{\k} (E',M_i)  \ar[r] 
&  \Hom_{\k} (E',M_{i+1}) \ar[r]& \cdots
}
\end{equation}

In Theorem 4.3 of \cite{EFS}, Eisenbud, Fl{\o}ystad and Schreyer show
that if $M$ is a graded $S$-module, with linear free resolution given
by $\mathbf{L}(P)$, then the dimension of $\Tor^{E'}_i(P,\k)$ can
be computed from the dimensions of the graded pieces of the local
cohomology modules of $M$. This result will be the main tool in
our application of BGG; for completeness (and because of the grading
differences) we prove a variant of their result in \S \ref{subs:torandext} 
below.

\subsection{The Eisenbud-Popescu-Yuzvinsky resolution}
\label{subs:epy}
Let $\A=\{H_1,\dots ,H_n\}$ be a central arrangement in $\C^{\ell}$, 
with complement $X(\A)$. Identify the cohomology ring 
$H^*(X(\A),\k)$ with the Orlik-Solomon algebra $A=E/I$, 
where $E$ is the exterior algebra on $V=\k^n$. The BGG correspondence 
was used in this context by Eisenbud-Popescu-Yuzvinsky \cite{EPY} 
to establish that $H_*(X(\A),\k) \cong \ann(I)$ has a linear free 
resolution over $E$; results on the first differentials in this resolution 
appear in \cite{DY}.

Fix a basis $e_1,\dots ,e_n$ for $V$, and let $x_1,\dots,x_n$ 
be the dual basis for $V^*$.   From \cite{EPY}, Corollary 3.2  
we have an exact sequence of $S$-modules:
\begin{equation}
\label{fa}
\xymatrixcolsep{21pt}
\xymatrix{
0\ar[r] & A_0  \otimes S \ar[r]^{d_0} &
A_1\otimes S \ar[r]^(.54){d_1} 
& \cdots \ar[r]^(.42){d_{\ell-1}}  & 
A_{\ell} \otimes S \ar[r] & F(A) \ar[r] & 0
}.
\end{equation}
The key point here is that the {\it complex} obtained by applying
BGG to $A$ is in fact exact, hence a free resolution of $F(A)$.
Notice that the differential 
$d_i\colon p\otimes 1 \mapsto \sum_{i=1}^n e_i  p \otimes x_i $
is precisely the differential in the Aomoto complex 
$(A,a)$, where the maps are given by multiplication by a 
generic linear form of the exterior algebra; the $x_i$ are 
simply the coefficients of this form.
With grading convention that $E$ is generated in degree one, 
$F(A)$ is generated in degree $-\ell$. With this choice, 
$\Ext_S^{\ell-1}(F(A),S)$ is generated in degree $\ell-1$; 
we will see in a bit that this is consistent
with the topological formula \eqref{eq:linmas}. 

\subsection{Relating $\Tor^E$ and $\Ext_S$}
\label{subs:torandext}
The following lemma is a restatement of Theorem 4.3 of \cite{EFS}. 
For our purposes it is preferable to use $\Ext^i(\bullet,S)$ rather than
the local cohomology modules $H^i_\m(\bullet)$; by local duality 
(\cite{E}, Theorem A.4.2) these modules encode essentially 
the same information.
\begin{lemma}
\label{mij}
With notation as above, we have:
\[
\dim_{\k} \Tor^E_i(A,\k)_j = \dim_{\k} \Ext_S^{i+\ell-j}(F(A),S)_{j}.
\]
\end{lemma}
\begin{proof}
As before, write $A'$ for the Orlik-Solomon algebra, with 
unit in degree $\ell$, and generators in degree $\ell-1$. 
A straightforward translation shows that 
\begin{equation}
\label{eq0}
\dim_{\k} \Tor^E_i(A,\k)_j = \dim_{\k} \Tor^{E'}_i(A',\k)_{\ell-j}
\end{equation}

Let $C$ denote the Cartan complex; $C$ is the minimal free 
resolution of $\k$ as an $E'$-module. For any $E'$-module $P$,  
\begin{equation}
\label{eq1}
\Tor^{E'}_i(P,\k) = H_i(P\otimes C).
\end{equation}
By Proposition 7.9 and Exercise 7.7 of \cite{e03},
$H_k(\Hom_S(\mathbf{L}(P),S))_{i+k}$ is 
dual to $H_i(P\otimes C)_{-i-k}$. If we apply 
this with $P = A'$, then \cite{EPY} tells us that
$\mathbf{L}(A')$ is a free resolution of $F(A')$, 
so that (\cite{E}, A3.11)
\begin{equation}
\label{eq2}
H_k(\Hom_S(\mathbf{L}(A'),S)) = \Ext^{-k}(F(A'),S).
\end{equation}

The reason for tensoring so that $A'$ has no components of negative
degree is that it keeps the indexing simple; in particular, we find that 
\begin{equation}
\label{eq3}
\dim_{\k} \Tor^{E'}_i(A',\k)_{\ell-j} = \dim_{\k} \Ext^{i+\ell-j}_S(F(A'),S)_{j-\ell}.
\end{equation}

Finally, since $F(A')$ is generated in degree zero, $F(A) =
F(A')\otimes S(\ell)$ is generated in degree $-\ell$, as claimed. So
\begin{align}
\label{eq4}
\dim_{\k} \Hom_S(F(A),S)_t 
&= \dim_{\k} \Hom_S(F(A')\otimes S(\ell),S)_t\\
&= \dim_{\k} \Hom_S(F(A'),S)_{t-\ell}. \notag
\end{align}
Thus, 
\begin{equation}
\label{eq5}
\dim_{\k} \Ext^{i+\ell-j}_S(F(A'),S)_{j-\ell} = \dim_{\k}
\Ext^{i+\ell-j}_S(F(A),S)_{j}. 
\end{equation}
Putting together \eqref{eq0}, \eqref{eq3}, and \eqref{eq5} 
finishes the proof.
\end{proof}

Recall that the regularity of an $E$-module $M$ is the smallest
integer $n$ such that $\Tor^E_i(M,\k)_{j} = 0$ 
for all $i,j$ with $j \ge i+ n+1$.
Lemma 2.3 of \cite{SS} is that $A$ is $(\ell-1)$-regular; 
this also follows from Lemma \ref{mij} and the 
EPY resolution of $F(A)$. For example, if we write 
$e_{ij}$ for $\dim_{\k} \Ext_S^{i}(F(A),S)_{j}$, 
then the minimal free resolution of $A$ over $E$ (for $\ell=3$) is:
\[
\xymatrixcolsep{15pt}
\xymatrix{
0& A \ar[l]  & E \ar[l] & E^{e_{22}}(-2) \oplus E^{e_{13}}(-3)\ar[l] 
&E^{e_{23}}(-3) \oplus  E^{e_{14}}(-4)\ar[l] 
& \cdots \ar[l]
}
\]

\subsection{The module $\B$ as an $\Ext$ module}
\label{subs:bext}
We now connect the players.  Recall that the  linearized Alexander 
invariant (over $\k$) is the $S$-module defined via the exact sequence
\begin{equation*}
\label{eq:BB}
\xymatrix{
(A_2 \oplus E_3)\otimes S \ar[rr]^(.55){\alpha_2 \otimes \id +\delta_3}
&& E_2\otimes S \ar[r] &  \B  \ar[r] & 0
},
\end{equation*}
where $\alpha_k\colon A_k\to E_k$ is the adjoint of the canonical 
surjection $\gamma_k\colon E_k\surj A_k$ (recall, we write 
$\B$ instead of $\B\otimes \k$, for simplicity).   

\begin{thm}
\label{bext2}
The linearized Alexander invariant is functorially determined by the 
Orlik-Solomon algebra, as follows:
\begin{equation}
\B \cong\Ext_S^{\ell-1}(F(A), S).
\end{equation}
\end{thm}

\begin{proof}
Consider the following commuting diagram, with exact rows and columns:
\begin{equation}
\label{eq:bigcd}
\xymatrixcolsep{20pt}
\xymatrix{
&& & 0 \ar[d] &  & 0 \ar[d] & 0 \ar@{.>}[d]  \\
&&0 \ar[r]  & I_2\otimes S \ar[r]^{\partial_3^*} \ar[d]& \cdots 
\ar[r]^{\partial_{\ell-1}^*} & I_{\ell}\otimes S \ar[r]
\ar[d] &F(I) \ar[r] \ar@{.>}[d] & 0 \\
0 \ar[r] &E_0\otimes S  \ar[r]^{\delta_1^*} \ar[d]^{=} 
&E_1\otimes S \ar[r]^{\delta_2^*} \ar[d]^{=}
& E_2\otimes S \ar[r]^{\delta_{3}^*} \ar[d]^{\gamma_2} 
& \cdots \ar[r]^(.42){\delta_{\ell-1}^*}  & 
E_{\ell} \otimes S \ar[r] \ar[d]^{\gamma_{\ell}} &F(E) 
\ar[r] \ar@{.>}[d] & 0\\
0 \ar[r] & A_0\otimes S \ar[r]^{d_0} &A_1\otimes S \ar[r]^{d_1} & 
A_2\otimes S\ar[d] \ar[r]^{d_2}  & \cdots \ar[r]^(.42){d_{\ell-1}}  & 
A_{\ell} \otimes S \ar[d] \ar[r] & F(A) \ar[r] \ar@{.>}[d] & 0\\
&&& 0 & & 0& 0}
\end{equation}

The rows in the diagram are EPY resolutions.  Notice that the middle 
row is just a truncation of the Koszul complex over $S$.  
All columns but the last one (with solid arrows) are exact by 
definition of the Orlik-Solomon algebra $A=E/I$. The exactness of
the two bottom rows, combined with the long exact sequence in homology,
shows that the rightmost column (marked with dotted arrows) is
exact. This short exact sequence yields a long exact sequence of 
$\Ext$ modules:
\[
\xymatrixcolsep{16pt}
\xymatrix{
\cdots \ar[r] & \Ext_S^i(F(E),S) \ar[r] 
&    \Ext_S^i(F(I),S) \ar[r] & \Ext_S^{i+1}(F(A),S) \ar[r] 
&  \cdots
 }
\]
The Koszul complex is exact and self-dual, so 
$\Ext_S^i(F(E),S)=0$, for all $i<\ell$, thus 
$ \Ext_S^{i+1}(F(A),S) \cong \Ext_S^i(F(I),S)$, for $i<\ell-1$.    
Since all the vertical exact sequences (except the last one) 
consist of free $S$-modules, the dual sequences are also exact, 
and so $E_k/(\im \alpha_k)\cong (I_k \otimes S)^*$.  Hence:
\begin{align*}
\B&=E_2\otimes S/\im(\alpha_2\otimes \id +\delta_3) \\
&\cong(I_2 \otimes S)^*/\im(\partial_3)\\
&\cong \Ext_S^{\ell-2}(F(I),S)\\
&\cong \Ext_S^{\ell-1}(F(A),S).  
\qquad\qquad \qed
\end{align*}
\renewcommand{\qed}{}\end{proof}

If char $\k=0$ (which recall is our standing assumption), then
we obtain:
\begin{cor}
\label{thetator}
The Chen ranks of an arrangement group $G$ equal the graded Betti 
numbers of the linear strand of the cohomology ring $A$ over the 
exterior algebra $E$:
\begin{equation}
\label{chentor}
\theta_k(G) = \dim_{\k} \Tor_{k-1}^E(A,{\k})_{k}, \quad \text{for all $k\ge 2$}.
\end{equation}
\end{cor}
\begin{proof}
We have:
\begin{align*}
\theta_k(G) 
&=\dim_{\k} \B_k &\qquad \text{by Theorem  \ref{fact2}} \\
&=\dim_{\k} \Ext_S^{\ell-1}(F(A),S)_k &\qquad \text{by Theorem \ref{bext2}} \\
&=\dim_{\k} \Tor^E_{k-1}(A,\k)_k &\qquad \text{by Lemma \ref{mij}}.
\end{align*}
\end{proof}

In particular, this shows that Conjectures \ref{conj:chen}  and 
\ref{conj:linconj} and are equivalent. We note that it is possible to assemble 
a different proof of Corollary \ref{thetator} using recent results of 
Fr\"oberg and L\"ofwall \cite{FL} on Koszul homology 
and homotopy Lie algebras.

\subsection{Chen groups of graphic arrangements}
\label{subs:chengraph}
We conclude this section with an application to a particularly nice 
class of arrangements. Given a simple graph $\Gamma$, with 
vertex set $\mathcal{V}=\{1,\dots , n\}$ and edge set $\mathcal{E}$, 
the corresponding {\em graphic arrangement}, 
$\A(\Gamma)=\{H_e\}_{e\in \mathcal{E}}$, consists of the hyperplanes 
in $\C^n$ of the form $H_{e}=\{z \in \C^n \mid z_i -z_j =0\}$ for which 
$e=(i,j)$ belongs to $\mathcal{E}$.  

For example, if $\Gamma=K_n$ is the complete graph on $n$ vertices, 
then $\A(K_n)$ is the braid arrangement in $\C^n$.  The resonance 
components of $\A(K_n)$ are in one-to-one 
correspondence to the $3$-vertex and $4$-vertex subsets of 
$\mathcal{V}$, and are all $1$-dimensional; see \cite[\S6.8]{CScv}, 
and compare with Example \ref{K4}.  

Now, since every graph on $n$ vertices  is a subgraph of the 
complete graph $K_n$, every graphic arrangement is a 
sub-arrangement of the braid arrangement $\A(K_n)$.  
It follows that $\RR^1(\A(\Gamma))$ has precisely 
one, $1$-dimensional component for each triangle 
or complete quadrangle in $\Gamma$.  

Combining Corollary~\ref{thetator} with our previous results 
from \cite{SS}, where the linear strand of a graphic arrangement 
is determined, we see that the Chen ranks conjecture holds for 
graphic arrangements.  More explicitly, we have the following.  
  
\begin{thm}
\label{graphic}
Let $\Gamma$ be a graph, $\A=\A(\Gamma)$ the corresponding 
graphic arrangement, and $G=\pi_1(X(\A))$ its group.  Then
\begin{equation}
\label{chengraph}
\theta_k(G) = (k-1)  \left( \kappa_2 + \kappa_3 \right),\quad  
\text{for all $k\ge 3$}, 
\end{equation}
where $\kappa_s$ denotes the number of complete subgraphs 
on $s+1$ vertices.  
\end{thm}

\begin{example}
\label{chenpn}
If $\A=\A(K_n)$ is the braid arrangement in $\C^n$, then 
$G$ is isomorphic to $P_n$, the pure braid group on $n$ strings. 
Applying formula \eqref{chengraph}, we find:
\begin{align*}
\theta_k(P_n) & = (k-1) \cdot \left( \binom{n}{3} + \binom{n}{4} \right) \\
& = (k-1) \cdot \binom{n+1}{4}, \qquad \text{for all $k\ge 3$}. 
\end{align*}
The Chen ranks of the pure braid groups were first computed  
in \cite{CS1}, using an arduous Gr\"obner basis 
computation. The above computation recovers the result of \cite{CS1}.
 \end{example}

\section{The rate of growth of the Chen groups}
\label{commalg}
It follows from Theorem \ref{fact2} that for $k \gg 0$ the Chen 
rank $\theta_k$ is given by a rational polynomial $P(k)$. In this
section, we prove that the degree of $P(k)$ is equal to the dimension 
of $\RR^1(\A)$.   As an application, we compute the complexity 
of the Orlik-Solomon algebra $A$, viewed as a module over the exterior 
algebra $E$, in the case when $\A$ is a central arrangement  in $\C^3$. 

\subsection{Primary decomposition of $\ann(\B)$}
\label{subs:prim}
We know that the Fitting ideal of the presentation matrix 
of $\B$ has the same radical as the annihilator
of $\B$, and by \cite{CScv}, the corresponding 
variety is precisely $\RR^1(\A)$. This variety has a number of nice 
properties:  as noted earlier, Libgober and Yuzvinsky \cite{LY} 
showed that it is the union of disjoint projective subspaces, 
each of dimension at least one.  

For each component $L\subset \RR^1(\A)$, denote by 
$\p=\p_{L}$ the corresponding minimal prime ideal in $S$.  
Slightly abusing notation, we write $\PP$ for the set of such ideals.   
For each $\p \in \PP$, denote by  $L=L_{\p}$ the 
corresponding linear subspace of $\RR^1(\A)$.  
Then:
 \begin{equation}
\label{annb}
\RR^1(\A)=\bigcup_{\p\in \PP}L_{\p} \qquad \text{and} \qquad 
\sqrt{\ann(\B)} = \bigcap_{\p\in \PP}  \p.
\end{equation}

For each $\p\in \PP$, there is at least one homogeneous element 
$m \in \B$ whose annihilator is $\p$, i.e., $\langle m\rangle  \cong S/\p$, 
as $S$-modules. Let $M(\p)$ denote the submodule of $\B$ 
annihilated by $\p$:
\begin{equation}
\label{eq:mp}
M(\p) = \{x \in \B\mid \ann (x)= \p\}.  
\end{equation}
It is easy to see that there is almost no interaction between 
these (finitely-generated) submodules. 
\begin{lemma}
\label{lem:fgmod}
Let $N$ be a finitely-generated, graded $S$-module, and 
$n_1, n_2 \in N$ satisfy $n_1 \notin \langle n_2\rangle$ and 
$n_2 \notin \langle n_1\rangle$. Suppose
$\ann \langle n_i\rangle = P_i  \subsetneq \m $
and $P_1 + P_2 = \m$, where $\m$ is
the maximal homogeneous ideal. Then
\begin{equation}
\langle n_1\rangle \cap  \langle n_2\rangle = 0.
\end{equation}
\end{lemma}
\begin{proof}
Suppose $x = a_1 n_1 = a_2 n_2$, with $a_i \in S$.
Then obviously both $P_1$ and $P_2$ annihilate $x$, so the
maximal ideal $\m$ kills $x$. By degree considerations,
we must have $a_i\in \m$.
Since $\m\cdot a_1 n_1 =0$, we find that
$a_1^2 n_1 = 0$.  But $P_1$ is prime, so
$a_1^2 \in P_1$ implies $a_1 \in P_1$.  
Hence $x = 0$.
\end{proof}

\subsection{Local cohomology and regularity}
\label{subs:loccoho}
Associated to any finitely generated, graded module $N$, 
there is an exact sequence (see Eisenbud \cite{E}):
\begin{equation}
\label{eq:exact}
\xymatrixcolsep{18pt}
\xymatrix{
0  \ar[r] & H^0_\m(N)   \ar[r] & N
  \ar[r] &\sum_{v \in {\Z}}H^0(\widetilde N(v))    \ar[r] &
H^1_\m(N)   \ar[r] & 0
},
\end{equation}
where $H^*_{\m}(N)$ denotes the local cohomology of $N$, 
supported at the maximal ideal $\m$, and $\widetilde N$ 
is the sheaf associated to $N$.

In \cite{Yuz}, Yuzvinsky showed that the module $F(A)$ is annihilated by 
$f=\sum_{i=1}^{n} x_i$; thus $F(A)_f=0$.  Localization commutes 
with the $\Ext$ functor, so for any $i$, $\Ext_S^i(F(A),S)_f=0$, 
and we conclude that $f$ annihilates $\B= \Ext_S^{\ell-1}(F(A),S)$. 
In particular every element of $\B$ is torsion,
hence every element of $\B$ belongs to some $M(\p)$.
Combined with Lemma \ref{lem:fgmod}, this shows that the 
sheaf associated to $\B$ decomposes as a direct
sum of the sheaves associated to the modules $M(\p)$: 
\begin{equation}
\label{eq:dec}
\widetilde{\B} \cong \bigoplus_{\p\in\PP} \widetilde{M}(\p).   
\end{equation}

The modules $H^0_\m(\B)$ and $H^1_\m(\B)$ vanish in 
high degree (see Chapter 9 of \cite{e03}), so the exact 
sequence above yields a coarse approximation to the 
Chen ranks conjecture: for $k > \reg(\B)$,
\begin{equation}
\label{chenhp}
\theta_k (G)= \HP(\B,k) = \sum_{p\in \PP} \HP(M(\p),k).
\end{equation}

\begin{lemma}
\label{mp}
Let $\p$ be a prime ideal in $\PP$.  Then 
$\deg \HP(M(\p),k)= \dim L_{\p}$.  
\end{lemma}

\begin{proof}
The module $M(\p)$ is generated by finitely many homogeneous 
elements, each having annihilator $\p$.  Let $m$ be one of those 
elements, and $\langle m \rangle$ the submodule it generates. 
 Then  $ \langle m \rangle \cong S(-d)/\p $, where $d$ 
is the degree of $m$. On the other hand, $L_{\p}=V(\p)$ is 
a linear subspace of $E_1$, say of dimension $r$. Hence:
\begin{equation}
\label{hpm}
\HP(\langle m\rangle,k)  = \HP(S(-d)/\p,k)
=\binom{r +k -d}{k - d} =\frac{k^{r}}{r!}+ \cdots ,
\end{equation}
and so $\deg \HP(\langle m\rangle,k)=r$.  
Since $\langle m \rangle$ is a submodule of $M(\p)$, 
this implies  
\[
\deg \HP(M(\p),k)\ge r.
\]
On the other hand, $M(\p)$ is generated by finitely many such 
submodules, and so \eqref{hpm} also implies $\deg \HP(M(\p),k)\le r$. 
\end{proof}

Combining Theorem \ref{fact2}, Equation \eqref{eq:dec}, and 
Lemma \ref{mp}, we obtain the following.

\begin{thm}
\label{chenpoly}
There exists a polynomial $\HP(t)\in \Q[t]$, of degree 
equal to the dimension of the projective variety $\RR^1(\A)$, such that 
\begin{equation}
\theta_k(G) = \HP(k),\quad \text{for all $k > \reg(\B)$}. 
\end{equation}
\end{thm}

\subsection{The complexity of the OS-algebra}
\label{subs:cxos}
Let $M$ be a finitely generated module over a local or graded ring $R$ 
with residue field $\k$, and let 
\begin{equation}
\label{torbn}
\beta_i(M)=\dim_\k \Tor^R_i(M,\k)
\end{equation}
be the rank of the $i^{\text{th}}$ free module in a minimal free resolution of
$M$ over $R$. We then have the following measure of the 
growth of the numbers $\beta_i(M)$, see \cite{A}. 

\begin{definition}
\label{def:complexity}
The {\em complexity} of the module $M$ is the integer 
\[
\cx_R(M) = \inf \{d \in \Z \mid 
\exists\, f(t) \in \R[t] \text{ of degree $d-1$, with 
$\beta_i(M) \le f(i),\, \forall i \ge 1$} \}.
\]
\end{definition}

Now let $\A$ be a central arrangement of $n$ planes in $\C^3$.  
We wish to compute the complexity of the Orlik-Solomon algebra 
$A=E/I$, viewed as a module over $E$.  We start with an 
observation about the resonance variety $\RR^1(\A)$.   

\begin{lemma}
\label{lem:pencil}
If $\A$ is not a pencil, then $\dim \RR^1(\A) < n-2$.
\end{lemma}

\begin{proof}
If $\dim \RR^1(\A) \ge n-2$, then since $\RR^1(\A)$ is the  
union of a subspace arrangement in $\P^{n-1}$, it must 
contain a hyperplane $H\cong \P^{n-2}$. 
 But since the components of $\RR^1(\A)$ are linear
subspaces of dimension at least one, B\'{e}zout's theorem says 
any component $L $ of $\RR^1(\A)$ different from $H$ 
must meet $H$, contradicting the fact that $\RR^1(\A)$ 
consists of disjoint subspaces (notice that applying 
B\'{e}zout makes sense, because linear subspaces 
are disjoint over $\k$ iff they are disjoint over $\bar \k$).
Thus $\RR^1(\A) = H$, and so $\RR^1(\A)$ contains a 
single component (necessarily local), of dimension $n-2$.  
Hence $\A$ is a pencil.
\end{proof}

Write $b'_{i,j}=\dim_{\k} \Tor^E_i(A,\k)_j$.  We noted in \cite{SS} 
that the regularity of $A$ as an $E$-module is two; 
so $\beta_i(A) = b'_{i,i+1}+b'_{i,i+2}$.
Using the (minimal, free) resolution of $A$ over $E$, 
we can then write the Hilbert series of $A$ as 
\[
\sum_{i=0}^3 b_it^i =  (1+t)^n \sum_{i \ge 0}(-1)^i
(b'_{i,i+1}t^{i+1}+b'_{i,i+2}t^{i+2}).
\]
Since $b_0=1$, $b_1=n$, and $b_3=b_2-n+1$ (see \cite{OT}), 
we have
\begin{equation}
\label{hilbetti}
\frac{1+nt+b_2t^2+(b_2-n+1)t^3}{(1+t)^n} = \sum_{i \ge 0}(-1)^i
(b'_{i,i+1}t^{i+1}+b'_{i,i+2}t^{i+2}).
\end{equation}
Expanding the left hand side of \eqref{hilbetti}, we find it equals:
\[
1+\sum_{m=2}^{\infty}(-1)^{m} \binom{n+m-4}{m-2} \cdot
\frac{3 (m-2)+m b_2-(2 m-5) n-n^2}{m} t^m.
\]
Assume $m$ is odd, and expand the coefficient of $t^m$ as a polynomial
in $m$. We find that it is equal to 
$\frac{-b_2+2n-3}{(n-2)!}m^{n-2}+\cdots$. 
On the other hand, the coefficient of $t^m$ on the right hand side 
of \eqref{hilbetti} is  $b'_{m-1,m}-b'_{m-2,m}$.   
Now recall from \eqref{chentor} that $\theta_k=b'_{k-1,k}$, 
for all $k\ge 2$. Thus, we obtain
\begin{equation}
b'_{m-2,m}=\theta_{m-1}+ \frac{b_2-2n+3}{(n-2)!}m^{n-2}+\cdots
\end{equation}
A similar argument applies when $m$ is even. 
Notice that $b_2-2n+3=0$ iff $\A$ is a near pencil, i.e., $n-1$ planes
through a line, and one additional plane in general position. 

\begin{prop}
\label{prop:cx}
If $\A$ is a central arrangement of $n$ planes in $\C^{3}$, and 
$\A$ is not a near pencil, then $\cx_E(A) = n-1$. If $\A$ is a near
pencil, then $\cx_E(A) = n-2$.
\end{prop}

\begin{proof}
If $\A$ is a near pencil, then the resolution of $A$ over $E$ 
is linear, and we know that $\beta_m(A) = \theta_m$ is a 
polynomial in $m$ of degree $n-3$. 
If $\A$ is a pencil, then the resolution is also linear, and
$\beta_m(A) = \theta_m$ is a polynomial of degree $n-2$. If $\A$ is
not a pencil or near pencil, then by Lemma \ref{lem:pencil}, 
$\dim \RR^1(\A) < n-2$. It follows from Theorem \ref{chenpoly} 
that $\theta_m$ is a polynomial of degree at most $n-3$, 
so $\beta_m(A)$ is a polynomial of degree $n-2$, for $m\gg 0$. 
\end{proof}

\begin{remark}
\label{rem:cx}
In Theorem 3.1 of \cite{AAH}, Aramova, Avramov, and Herzog show 
that the complexity of an $E$-module is equal to the dimension of 
the (affine) rank variety; in Corollary 2.3 of \cite{EPY}, Eisenbud, 
Popescu, and Yuzvinsky prove that when the module in question 
is the Orlik-Solomon algebra, the codimension of the rank variety 
is equal to the number of central factors in an irreducible 
decomposition of $\A$.  On the other hand, near pencils are the 
only reducible arrangements in rank $\ell=3$.  So Proposition 
\ref{prop:cx} is in fact an instance of a more general result. 
\end{remark}

\section{A lower bound for the Chen ranks}
\label{sect:exseq}
In this section we prove a lower bound on the Chen ranks of 
an arrangement group. We continue to use the notation of the 
previous section.   

The idea of the proof is as follows. First, we choose an irreducible 
component of $\RR^1(\A)$ (we do not distinguish between the 
local and non-local components), which corresponds to a 
sub-ideal of $I_2$, generated by decomposable elements. 
A short exact sequence relates this sub-ideal 
of decomposable elements to $I$. After appealing to the 
Bernstein-Gelfand-Gelfand correspondence and dualizing, 
we obtain a long exact sequence of $\Ext$-modules. Finally, 
knowledge of the geometry of $\RR^1(\A)$, combined 
with a localization argument, yields the bound. 

\subsection{The modules $\B(\p)$}
\label{subs:bpmod}
Let $\p \in \PP$ be a minimal prime ideal, and 
$L_{\p}=V(\p) \subset E_1$ the corresponding 
linear subspace.  As shown by Falk in \cite{Fa},  Corollary 3.11, 
the following holds:
If $a,b \in L_{\p}$, then $a\wedge b \in I_2$. Thus, if we define  
\begin{equation}
\label{eq:ip}
I(\p) =\langle a \wedge b \mid a, b \in L_{\p} \rangle
\end{equation}
to be the ideal of $E$ generated by wedge products of 
pairs of elements from $L_{\p}$, then $I(\p)$ is a sub-ideal of $I$.  

Recall from Theorem \ref{bext2} that 
$\B \cong\Ext_S^{\ell-1}(F(A), S)$.
By analogy, define the {\em linearized Alexander 
invariant at $\p$} to be the $S$-module 
\begin{equation}
\label{eq:bp}
\B(\p) = \Ext_S^{\ell-2}(F(I(\p)),S). 
\end{equation}

\begin{example}
\label{ex:local}
Let $\A=\{H_0,\dots ,H_n\}$ be a pencil of lines in $\C^2$. 
Let $E$ be the exterior algebra on $e_0,\dots ,e_n$.  
The Orlik-Solomon ideal $I$ is generated by  
$\{ \partial e_{0ij} \}$, 
while the resonance variety  $\RR^1(\A)$ consists of a 
single component $L_{\p}= \spn\{ e_0-e_i\}$, 
with $\p= \langle \sum_{i=0}^n x_i, x_{r+1}, \ldots, x_{n} \rangle$   
and $I({\p})=I$.  After a linear change of variables in $E$, 
the ideal $I$ corresponds to $I_0=\langle e_1, \dots , e_n\rangle^2$,  
while the linear subspace $L_{\p}$ corresponds 
to  $L_0=\spn\{ e_1, \dots ,e_n\}$. 
A standard calculation (compare \cite{CSai}) shows 
that the module $\B(\p)$ has Hilbert polynomial 
$\HP(\B(\p),k)=(k-1) \binom{n +k-2}{k}$. 
\end{example}

\begin{example}
\label{ex:braidres}
Let $\A$ be the braid arrangement, discussed in Example 2.3.
The generators of $I_2$ are $f_1=\partial e_{145}$,
$f_2=\partial e_{235}$, $f_3=\partial e_{034}$, and
$f_4=\partial e_{012}$.   These give rise to various 
``local" data.  For example, 
$\p_1=\langle x_0,x_2,x_3,\sum x_i\rangle$,
$L_{\p_1}=\spn\{e_1-e_4,e_1-e_5\}$, 
$I(\p_1)=\langle (e_1 -e_4)(e_1-e_5)\rangle$, 
and $\B(\p_1) =\coker \left( \vartheta_{\p_1}\colon S^{4} \to S \right)$, 
where $\vartheta_{\p_1} =\begin{pmatrix} 
x_0 &  x_2  & x_3  & \sum x_i  \end{pmatrix}$. 

Now recall $\A$ supports a non-trivial neighborly partition, 
 $\Pi=(05|13|24)$.  The corresponding non-local
component $L_{\Pi} \subset \RR^1(\A)$ is spanned by 
$\eta_1=e_0-e_1-e_3+e_5$ and
$\eta_2=e_0-e_2-e_4+e_5$. 
The associated prime is 
$\p=\langle x_0-x_5, x_1-x_3,x_2-x_4, \sum x_i\rangle$, 
whereas 
$I(\p)=\langle \eta_1\wedge \eta_2\rangle$.  Finally,
$\B(\p) =\coker \left( \vartheta_\p \colon S^{4} \to S \right)$, where 
$\vartheta_{\p} =\begin{pmatrix} 
x_0-x_5 &  x_1-x_3  & x_2-x_4  & \sum x_i  \end{pmatrix}$. 
\end{example}

\begin{lemma}
\label{lem:ip}
If $\dim L_{\p}=\dim L_{\q}$, there is a linear automorphism of $E$  
taking $I(\p)$ to $I(\q)$. 
\end{lemma}

\begin{proof}  
After a suitable linear change of variables, we can assume 
$L_{\p}$ is the subspace spanned by $\{e_1, \dots ,e_r\}$, 
where $r=\dim L_{\p}$.  
Let $I(\p)$ be the corresponding ideal of $E$. 
We know that every element $e_i\wedge e_j$ ($1\le i<j\le r$) 
is in $I(\p)$.  Hence, $I(\p)=\langle e_1, \dots ,e_r\rangle ^2$.
\end{proof}

In particular, the sub-ideals $I(\p)$  (and thus, the modules $B(\p)$) 
do not distinguish between local and non-local resonance. 
In view of Example~\ref{ex:local}, we obtain the following. 

\begin{cor}
\label{cor:hpbp}
The Hilbert polynomial of $\B(\p)$ is given by:
\begin{equation}
\HP(\B(\p),k)=(k-1) \binom{\dim L_\p +k-1}{k}. 
\end{equation}
\end{cor}

\subsection{Snake lemma and localization}
\label{subs:snake}
For simplicity, we will assume from now on that 
$\A$ is a central arrangement in $\C^3$. 
For a fixed $\p\in \PP$, consider the short  exact 
sequence of $E$-modules
\[
\xymatrix{ 0 \ar[r] & I(\p) \ar[r] & I \ar[r] & I/I(\p) \ar[r] & 0 }.
\]
Using $I_1=0$ (which also implies $I(\p)_1=0$),
we obtain the following commuting diagram of $S$-modules
\begin{equation}
\label{eq:cd2}
\xymatrixcolsep{18pt}
\xymatrixrowsep{18pt}
\xymatrix{
& &0 \ar[d] & 0 \ar[d] &  & \\
  & 0  \ar[r]  & I(\p)_2\otimes S \ar[r] \ar[d]
& I(\p)_3\otimes S \ar[r] \ar[d] &F(I(\p)) \ar[r]   & 0 \\
   & 0  \ar[r]  & I_2\otimes S \ar[r] \ar[d]
& I_3\otimes S \ar[r] \ar[d] &F(I) \ar[r]   & 0 \\
0 \ar[r]  & K_1  \ar[r] 
& I_2/I(\p)_2\otimes S \ar[r] \ar[d]
& I_3/I(\p)_3\otimes S \ar[r] \ar[d] &F(I/I(\p)) \ar[r]  & 0 \\
& & 0 & 0 & }
\end{equation}
The Snake Lemma now yields the exact sequence
\[
\xymatrix{
0\ar[r] & K_1 \ar[r] & F(I(\p)) \ar[r] & F(I) \ar[r] & F(I/I(\p)) \ar[r] & 0
}.
\]

All the $S$-modules in the two middle columns of 
diagram \eqref{eq:cd2} are free. Thus, their duals 
are also free. Hence, dualizing the full
diagram, we obtain the commuting diagram
\begin{equation*}
\label{eq:cd3}
\xymatrixcolsep{11pt}
\xymatrixrowsep{16pt}
\xymatrix{
& &0  & 0  &  & \\
0 &\Ext_S^1(F(I(\p)),S) \ar[l]  & \left(I(\p)_2\otimes S\right)^* 
\ar[l] \ar[u]
& \left(I(\p)_3\otimes S\right)^* \ar[l] \ar[u] &F(I(\p))^* \ar[l]  & 
0\ar[l]  \\
0&\Ext_S^1(F(I),S) \ar[l]  & \left(I_2\otimes S\right)^* \ar[l] \ar[u]
& \left(I_3\otimes S\right)^* \ar[l] \ar[u] &F(I)^* \ar[l]  & 0\ar[l]  \\
0  & K_2 \ar[l]  & \left(I/I(\p)_2\otimes S\right)^* \ar[l] \ar[u]
& \left(I_3/I(\p)_3\otimes S\right)^* \ar[l] \ar[u]
&F(I/I(\p))^* \ar[l]& 0\ar[l]  \\
& & 0 \ar[u]& 0 \ar[u]& }
\end{equation*}
from which we derive, as before, an exact sequence
\[
\xymatrixcolsep{12pt}
\xymatrixrowsep{10pt}
\xymatrix{
0\ar[r] & F(I/I(\p))^* \ar[r] & F(I)^* \ar[r] & F(I(\p))^* \ar@{->}@/_/[dl] \\
&&K_2 \ar[r] & \Ext_S^1(F(I),S) \ar[r] 
& \Ext_S^1(F(I(\p)),S)\ar[r] & 0.
}
\]

Let $\epsilon(\p)\colon \Ext_S^1(F(I),S) \surj \Ext_S^1(F(I(\p)),S)$
be the surjective map on the right.  By construction, the only 
associated prime of  $B(\p)=\Ext_S^1(F(I(\p)),S)$ is $\p$.  
Put $C=\coker \left(\oplus_{\p} \epsilon(\p)\right)$.
We obtain the exact sequence
\begin{equation}
\label{eq:extc}
\xymatrixcolsep{18pt}
\xymatrix{ \Ext_S^1(F(I),S) \ar[rr]^(0.36){\oplus_{\p} \epsilon(\p)}
&& \bigoplus_{\p \in \PP}\Ext_S^1(F(I(\p)),S) \ar[r] & C \ar[r] &0.}
\end{equation}

\begin{prop}
\label{prop:vanish}
The $S$-module  $C=\coker \left(\oplus_{\p} \epsilon(\p)\right)$ 
is supported only at the maximal ideal $\m$.
\end{prop}

\begin{proof}
Localize the sequence \eqref{eq:extc} at a (minimal) prime 
$\q\in \PP$.  All the summands in the middle term disappear, 
except for $\B(\q)=\Ext_S^1(F(I(\q)),S)$.  This is because 
$\ch\k= 0$, and so the components of $\RR^1(\A)$ are 
projectively disjoint, cf.~\cite{LY}. This gives an exact sequence
\begin{equation}
\label{eq:extloc}
\xymatrix{ \Ext_S^1(F(I),S)_{\q} \ar[rr]^{\epsilon(\q)_{\q}}
&& \Ext_S^1(F(I(\q)),S)_{\q} \ar[r] & C_{\q} \ar[r] &0}.
\end{equation}
Since localization is an exact functor, the surjection $\epsilon(\q)$
localizes to a surjection $\epsilon(\q)_{\q}$.  Thus, $C_{\q}=0$.

If $\q$ were an embedded prime for $\RR^1(\A)$, different from $\m$,
then the same argument as above shows that $C_{\q}=0$.  
This finishes the proof.  
\end{proof}

\subsection{A sharp inequality for $\theta_k$}
\label{subs:chineq}
We are now in position to prove one direction of 
Conjecture \ref{conj:chen}.  Let $G$ be an arrangement group.

\begin{thm}
\label{halfchen}
For $k$ sufficiently large,
\[
\theta_k(G) \ge (k-1) \sum_{\p \in \PP}  \binom{\dim L_{\p}+k-1}{k}.
\]
\end{thm}

\begin{proof}
Let $K=\ker \left(\oplus_{\p} \epsilon(\p)\right)$.
We have an exact sequence
\begin{equation*}
\label{eq:extker}
\xymatrixcolsep{13pt}
\xymatrix{0 \ar[r] & K\ar[r] & \Ext_S^1(F(I),S) 
\ar[rr]^(0.36){\oplus_{\p} \epsilon(\p)}
&& \bigoplus_{\p \in \PP}\Ext_S^1(F(I(\p)),S) \ar[r] & C \ar[r] & 0}.
\end{equation*}
By Theorem \ref{fact2} and Proposition \ref{prop:vanish}, 
we see that for $k \gg 0$, 
\begin{equation}
\label{hpfinal}
\begin{aligned}
\theta_k(G) & = \HP(\B,k)\\
& = \sum_{\p \in \PP}\HP(\Ext_S^1(F(I(\p)),S),k)+  \HP(K,k).
\end{aligned}
\end{equation}
By Corollary~\ref{cor:hpbp}:
\begin{equation}
\label{hplocal}
\HP(\Ext_S^1(F(I(\p)),S),k) = (k-1) \binom{\dim L_\p +k-1}{k}.  
\end{equation}
Combining Equations \eqref{hpfinal} and  \eqref{hplocal} 
yields the desired inequality.
\end{proof}

Theorem \ref{halfchen} gives the lower bound predicted by
the Chen ranks conjecture. Proving the other direction of the inequality 
is equivalent to showing that the module $K$ is of finite
length. While localization techniques should still be useful here, the
problem seems more delicate.  In the next section, 
we will indicate another possible way to attack this problem.

\section{Discussion and Examples}
\label{sect:examples}

\subsection{An alternate approach}
\label{subsect:alternate}
Let $\A$ be an arrangement in $\C^3$, and 
let  $\p\in \PP$ be one of the minimal prime ideals associated to the 
resonance variety $\RR^1(\A)$.  Recall we defined in \eqref{eq:mp} 
and  \eqref{eq:bp} two $S$-modules:  
\begin{equation}
M(\p)=\{x \in \B \mid \ann (x) = \p\} \qquad \text{and} \qquad 
\B(\p) = \Ext^1_S(F(I(\p)), S).  
\end{equation}
Computations suggest that $M(\p)$ and $\B(\p)$ define the same sheaf.  
Proving that $\widetilde{M}(\p)=\widetilde{\B}(\p)$ would imply, via 
\eqref{eq:exact}, the Chen ranks conjecture. If we put 
\begin{equation}
\B'=\bigoplus\limits_{\p \in \PP} \B(\p),  
\end{equation}
then $\widetilde{M}(\p)=\widetilde{\B}(\p)$  would also imply that 
\begin{equation}
  \B' =\sum_{v \in {\Z}}H^0(\widetilde \B(v)).
\end{equation}

\subsection{Almost neighborly partitions}
\label{subsect:almost}
The decomposable arrangements studied in \cite{CSai} have 
only local resonance. For these arrangements, $\B \cong \B'$, and 
it follows easily that $H^0_\m(\B)$ and $H^1_\m(\B)$ vanish. 
We note that while $H^0_\m(\B)$ and $H^1_\m(\B)$ are not important 
for computing the asymptotic Chen ranks, the next examples indicate
that they do capture interesting combinatorial features of an arrangement. 

\begin{example}
\label{delmac}

\begin{figure}[ht]
\setlength{\unitlength}{0.7cm}
\begin{picture}(5,4)(0.5,-0.5)
\put(3,3){\line(1,-1){3}}
\put(3,3){\line(-1,-1){3}}
\put(3,3){\line(0,-1){3}}
\put(1.5,1.5){\line(1,0){3}}
\put(0,0){\line(1,0){6}}
\put(0,0){\line(2,1){3}}
\put(6,0){\line(-3,1){4.5}}
\multiput(0,0)(6,0){2}{\circle*{0.3}}
\multiput(1.5,1.5)(1.5,0){3}{\circle*{0.3}}
\multiput(3,3)(0,-3){2}{\circle*{0.3}}
\put(2.4,1.2){\circle*{0.3}}
\put(0,0){\makebox(-1,0){$1$}}
\put(3,0){\makebox(0,-1){$2$}}
\put(6,0){\makebox(1,0){$5$}}
\put(1.5,1.5){\makebox(-1,0){$7$}}
\put(4.5,1.5){\makebox(1,0){$3$}}
\put(3,3){\makebox(0,1){$0$}}
\put(3,1.5){\makebox(0.6,0.6){$4$}}
\put(2.4,1){\makebox(0,-0.55){$6$}}
\end{picture}
\caption{\textsf{The deleted MacLane matroid}}
\label{fig:delml}
\end{figure}

The deleted MacLane arrangement corresponds to a matroid on $8$
points, obtained by deleting a line from the $\mathtt{ML}_8$ matroid, 
see Figure \ref{fig:delml}.   The resonance variety
has $7$ local components, corresponding to the triple points.  
A computation shows:
\begin{center}
$\begin{array}{c|cccc}
k&  \dim_\k H^0_\m(\B)_k&  \dim_\k \B_k & \dim_\k \B'_k&
\dim_\k H^1_\m(\B)_k\\[3pt]
\hline\\[-6pt]
2   &    0      &       7      &     7  & 0\\
3   &    1      &     15      &   14  & 0\\
4   &    0      &      21     &   21  & 0
\end{array}$
\end{center}
Moreover, 
$\Hilb(\B)=\tfrac{7t^2}{(1-t)^2}+t^3=7t^2+15t^3+21 t^4+28 t^5 +\cdots$,
and so $\theta_k=7(k-1)$, for $k\ge 4$.
\end{example}

The interesting feature in the above example is the jump of $1$ 
in the dimension of $\B_3$ (when compared to that of $\B'_3$), 
coming from the non-vanishing of $H^0_\m(\B)_3$. 
This jump can be explained by the partition $(06|13|27|45)$, 
which is ``almost'' neighborly, though not neighborly. 

\begin{definition}
\label{def:ANP}
A partition $\Pi$ of $\A$ is {\em almost neighborly} if, for
any rank two flat $Y\in L_2(\A)$ with $\mu(Y)>1$ and any 
block $\pi$ of $\Pi$,
\begin{equation}
\label{eq:almost np}
\mu(Y) \le \abs{Y\cap \pi} \Longrightarrow Y\subseteq \pi. 
\end{equation}
\end{definition}

While such partitions cannot contribute to resonance (at least 
not in characteristic $0$), it seems that they manifest their 
presence in the free resolution of $A$ over $E$. 
Computations indicate that non-vanishing $H^0_\m(\B)$ 
is often linked to the existence of almost neighborly partitions. 
Do there exist examples where $H^0_\m(\B)$ is nonzero, but 
which have no almost neighborly partitions?

\subsection{Regularity and asymptotic range}
\label{subsect:asymptotic}
The fact that the module $\B$ may not be saturated means 
that the problem of determining for which $k$ the asymptotic 
values for $\theta_k$ given by Conjecture  \ref{conj:chen} 
are valid is rather subtle. In particular, the regularity of $\B$ 
is related to the homology of the linear strand of $A$ over $E$ 
(see \cite{EFS}, Corollary 2.4), and will be influenced by the 
non-vanishing of $H^0_\m(\B)$.

Regularity also influences  the values for which the Hilbert 
polynomial and Hilbert function agree (see \cite{V}, Corollary B.4.1).  
This is illustrated in the next example.  

\begin{example}
\label{ex:ceva3}
\begin{figure}[ht]
\setlength{\unitlength}{0.47cm}
\begin{picture}(3,8.4)(0.5,-1.8)
\multiput(0,0)(2,0){3}{\line(0,1){4}}
\multiput(0,0)(0,2){3}{\line(1,0){4}}
\multiput(0,0)(0,2){2}{\line(1,1){2}}
\multiput(2,0)(0,2){2}{\line(1,1){2}}
\multiput(0,4)(0,-2){2}{\line(1,-1){2}}
\multiput(2,4)(0,-2){2}{\line(1,-1){2}}
\multiput(0,0)(0,2){3}{\circle*{0.4}}
\multiput(2,0)(0,2){3}{\circle*{0.4}}
\multiput(4,0)(0,2){3}{\circle*{0.4}}
\qbezier(0,0)(-5,10)(2,4)
\qbezier(0,2)(-5,-4)(4,0)
\qbezier(0,4)(9,9)(4,2)
\qbezier(2,0)(8,-5)(4,4)
\put(0.2,4){\makebox(-0.4,1.1){{\small $1$}}}
\put(2.2,4){\makebox(-0.4,1.1){{\small $2$}}}
\put(4.2,4){\makebox(-0.4,1.1){{\small $3$}}}
\put(0,2){\makebox(-0.5,1.3){{\small $4$}}}
\put(2,2){\makebox(-0.5,1.3){{\small $5$}}}
\put(4,2){\makebox(-0.5,1.3){{\small $6$}}}
\put(0,0){\makebox(-0.4,-1.1){{\small $7$}}}
\put(2,0){\makebox(-0.4,-1.1){{\small $8$}}}
\put(4,0){\makebox(-0.4,-1.1){{\small $9$}}}
\end{picture}
\caption{\textsf{The  $\text{Ceva}(3)$ matroid}}
\label{fig:ceva}
\end{figure}

The $\text{Ceva}(3)$ arrangement (also known
as the monomial arrangement $\A(3,3,3)$), is a realization
of the affine plane over $\Z_3$, with defining polynomial
$Q=(x^3-y^3) (x^3-z^3)(y^3-z^3)$.  Its matroid is depicted 
in Figure \ref{fig:ceva}.  The resonance variety
has $12$ local components, corresponding to the triple points,
and $4$ essential, $1$-dimensional components, corresponding
to the neighborly partitions $(123 | 456 | 789)$,
$(147 | 258 | 369)$, $(159 | 267 | 348)$, $(168 | 249 | 357)$.
We compute:
\begin{center}
$\begin{array}{c|cccc}
k &  \dim_\k H^0_\m(\B)_k&  \dim_\k \B_k & \dim_\k \B'_k &
\dim_\k H^1_\m(\B)_k\\[3pt]
\hline\\[-6pt]
2   &    0     &     12      &     16  & 4\\
3   &    8     &     40      &     32  & 0\\
4   &    8     &     56      &     48  & 0\\
5   &    0     &     64      &     64  & 0\\
\end{array}$
\end{center}
Also, 
$\Hilb(\B)=\tfrac{16t^2}{(1-t)^{2}}+4t^2(2t^2+2t-1)=
12t^2+40 t^3+56 t^4+64t^5+80 t^6+\cdots$, 
and so $\theta_k=16 (k-1)$, for $k\ge 5$.
\end{example}

As noted in \S\ref{subs:reschen}, the original form 
of Conjecture  \ref{conj:chen}  assumed that
taking $k\ge 4$ would suffice to ensure equality. 
The non-vanishing of $H^0_\m(\B)_4$ in Example \ref{ex:ceva3} 
shows that a larger $k$ is needed, in general.  

\subsection*{Acknowledgements}  Computations using 
Macaulay~2 \cite{GS} were essential to our efforts. We are  
grateful to MSRI, where part of this work was performed.  
Thanks to the referee for a very careful reading of the 
manuscript, and for useful suggestions.

\bibliographystyle{amsalpha}

\end{document}